\documentclass[11pt]{article}
\usepackage{latexsym,amsmath,amsfonts}
\begin{document}
\begin{center}
{\LARGE\textbf{Emergence in Random Noisy Environments}}\\
\bigskip
\bigskip
Yilun Shang\footnote{Department of Mathematics, Shanghai Jiao Tong University, Shanghai 200240, CHINA. email: \texttt{shyl@sjtu.edu.cn}}\\
\end{center}

\bigskip

\begin{abstract}
We investigate the emergent behavior of four types of generic
dynamical systems under random environmental perturbations.
Sufficient conditions for nearly-emergence in various scenarios are
presented. Recent fundamental works of F. Cucker and S. Smale on the
construction and analysis of flocking models directly inspired our
present work.



\end{abstract}

\bigskip
\bigskip
\noindent{\Large\textbf{1. Introduction}}
\smallskip

The emergent behaviors of a large number of autonomous interacting
agents such as flocking/schooling/swarming/herding in
birds/fishes/bacterica/animals \cite{6,7,18}, multi-agent
cooperative coordination in mobile networks \cite{16,17} and
emergence of a common language in primitive societies
\cite{5,11,12,15} have been attracting great research attentions
since the last two decades from biologists, computer scientists,
physicists, sociologists, engineers and mathematicians.

Recently,  Cucker and Smale \cite{1} have proposed a remarkable
model aiming to exploring the flocking phenomenon and mathematical
analysis is performed to show the convergence results only depend on
some initial states of the population. This notable feature is in
contrast with the previous models (e.g. the so called Vicsek model
\cite{18,7}) where convergence relies on the global behavior of the
agents' trajectories (or on the neighborhood graphs of the
underlying dynamical systems), which are quite hard to verify in
general. The same authors \cite{3} extend the model later to a more
general setting beyond flocking. \cite{13} further develops a
hierarchical leadership architecture in the Cucker-Smale flocking
model. The work in \cite{2} focuses on a situation where uniform or
Gaussian noises are involved in the environments. A hydrodynamic
description and the mean-field limit of this very model are also
provided in \cite{14}.

The starting point of our present work is directly motivated by the
aforementioned series work. Primarily, we want to refine the
rudimental results (in the noisy environment) in \cite{2} and extend
them to more general scenarios such as those discussed by \cite{3};
and try to shed some light on the understanding of various emergence
behaviors observed in diverse natural, social and man-made complex
systems \cite{20,23}. To do so, we first introduce four types of
non-autonomous, nonlinear dynamical systems; two ( I(D) and II(D) )
for discrete time and two ( I(C) and II(C) ) for continuous time. In
each case, we provide a convergence analysis. Systems I(D) and I(C)
are adapted from \cite{3,1} and the underlying idea stems from the
birds flocking in a noisy environment. Whereas the original idea
behind systems II(D) and II(C) is the linguistic evolution with some
possible fluctuations in a primitive society. The random noises
considered here may reflect the change of the environment which is
usually unclear to the objects. Moreover, information interaction
among agents may be contaminated or corrupted by errors. Hence, it
becomes significant to analyze systems in the presence of random
noises. We mention that the systems tackled in this paper are
quintessential in the sense of reflecting some typical mechanisms
behind emergence (see Remark 1 in Section 2.1), but by no means
limited to flocking or language evolution since we will treat them
in a quite general manner with emphasis on the methodology. Concrete
examples will follow and illustrate the strength of our generic
frameworks. Some other related work about emergent behaviors under
random environmental perturbation can be found in e.g.
\cite{9,10,19,24} and references therein.

The rest of this paper is organized as follows. In Section 2, we
will study the discrete time models I(D), II(D) and the statement of
main results in this section appears at Section 2.2. Section 3 is
devoted to the continuous counterparts I(C), II(C) and see Section
3.2 for the statement of main results in the continuous case. We
then draw our conclusion and discuss future direction in Section 4.

\bigskip

\noindent{\Large\textbf{2. Discrete-time emergence}}
\smallskip

Let $k\in\mathbb{N}$. We assume the population under consideration
consists of $k$ agents throughout the paper.

\smallskip
\noindent{\large\textbf{2.1 Models setup (I(D), II(D))}}
\medskip

We shall first introduce the dynamical system I(D), which is
developed similarly with that considered in \cite{3}.

Suppose $X$ and $Y$ are two given inner product spaces whose
elements are denoted as $x$ and $y$, respectively. Let
$x(t)=(x_1(t),\cdots,x_k(t))\in X$ and $(y_1(t),\cdots,y_k(t))\in
Y^k$ represent two kinds of characteristics of the agents at time
instant $t$. Convergence of $x\in X$ (or $y\in Y^k$) is naturally
understood as entrywise convergence as $t$ approaches infinity. Let
$\triangle$ signify the diagonal of $Y^k$, that is,
$\triangle:=\{(y,\cdots,y)|\ y\in Y\}$. Denote
$\widetilde{Y}:=Y^k/\triangle$ and fix an inner product
$\langle\cdot,\cdot\rangle$ in $\widetilde{Y}$, which induces a norm
$||\cdot||$. (Here in the discrete case, we do not really need an
inner; what we want is $\widetilde{Y}$ should be a normed space. The
same remark applies to $X$ and $Y$.) Since $\widetilde{Y}$ is a
finite dimensional space,
$\hat{y}:=(y_1(t),\cdots,y_k(t))\rightarrow(y_0,\cdots,y_0)$ for
some $y_0\in Y$ if and only if $||\bar{\hat y}-\bar0||\rightarrow0$,
where $\bar a:=a+\triangle\in \widetilde{Y}$ for $a\in Y^k$. In what
follows, we denote norms in all different spaces as $||\cdot||$ with
some ambiguity, but the proper meaning will be clear in the context.

For $x\in X$, $y\in \widetilde{Y}$, consider the following dynamical
system:
$$
I(D):\qquad\Big\{\begin{array}{l} x(t+h)=x(t)+hJ(x(t),y(t))\\
y(t+h)=S(x(t))y(t)+hH(t)
\end{array}
$$
Here $h$ is the time step and we shall denote in the sequel
$x[t]:=x(th)$, $y[t]:=y(th)$ and $H[t]:=H(th)$ for brevity. Take
$t\in\mathbb{N}$ herein. We now explain the notations in system
I(D). Let $J:X\times\widetilde{Y}\rightarrow X$ be a Lipschitz or
$C^1$ operator satisfying, for some $C,\delta>0$, $0\le \gamma<1$,
that
\begin{equation}
||J(x,y)||\le C(1+||x||)^{\gamma}||y||^{\delta}\label{1}
\end{equation}
for all $x\in X$, $y\in\widetilde{Y}$. Let $S:X\rightarrow
\mathrm{End}(\widetilde{Y})$ be an operator satisfying, for some
$G>0$, $\beta\ge0$, that
\begin{equation}
||S(x)||\le 1-\frac{hG}{(1+||x||)^{\beta}}\label{2}
\end{equation}
for all $x\in X$. The operator norm in (\ref{2}) is defined as
$||S(x)||=\sup_{y\ne0 \atop
y\in\widetilde{Y}}\frac{||S(x)y||}{||y||}$. Let
$H:(\Omega,\mathcal{F},P)\rightarrow(\widetilde{Y},\mathcal{B}(\widetilde{Y}))$
be a random element. $(\Omega,\mathcal{F},P)$ is some probability
space and $\mathcal{B}(\widetilde{Y})$ is the Borel $\sigma$-algebra
on $\widetilde{Y}$. We assume $H[t]$ is independent and identically
distributed for different $t\in\mathbb{N}$. Notice that
$||H||:(\Omega,\mathcal{F},P)\rightarrow(\mathbb{R},\mathcal{B}(\mathbb{R}))$
is a random variable and let $F(x):=P(||H||\le x)$ for
$x\in\mathbb{R}$ be the distribution function of $||H||$.

Next, we present our dynamical system II(D) as follows. The spaces
$X,Y,\widetilde{Y}$ are defined as before. For $x\in X$, $y\in
\widetilde{Y}$, consider the following dynamical system:
$$
II(D):\qquad\Big\{\begin{array}{l} x(t_1+h_1)=S_1(y(t_2))x(t_1)+h_1H_1(t_1)\\
y(t_2+h_2)=S_2(x(t_1))y(t_2)+h_2H_2(t_2)
\end{array}
$$

Here $h_1,h_2$ are the time steps w.r.t. $x$ and $y$. Take $t_1,
t_2\in\mathbb{N}$ and we denote $x[t_1]=x(t_1h_1)$,
$y[t_2]=y(t_2h_2)$, $H_1[t_1]=H_1(t_1h_1)$ and
$H_2[t_2]=H_2(t_2h_2)$. In light of these notations, the system
II(D) can be rewritten as follows
$$
\qquad\Big\{\begin{array}{l} x[t+1]=S_1(y[t])x[t]+h_1H_1[t]\\
y[t+1]=S_2(x[t])y[t]+h_2H_2[t]
\end{array}
$$

In analogy with the system I(D), suppose
$S_1:\widetilde{Y}\rightarrow\mathrm{End}(X)$ and
$S_2:X\rightarrow\mathrm{End}(\widetilde{Y})$ are two operators
satisfying, for some $G_1,G_2>0$, $\beta_1, \beta_2\ge0$, that
\begin{equation}
||S_1(y)||\le 1-\frac{h_1G_1}{(1+||y||)^{\beta_1}},\qquad
||S_2(x)||\le 1-\frac{h_2G_2}{(1+||x||)^{\beta_2}}\label{3}
\end{equation}
for all $y\in\widetilde{Y}$, $x\in X$. Let
$H_1:(\Omega_1,\mathcal{F}_1,P_1)\rightarrow(X,\mathcal{B}(X))$ and
$H_2:(\Omega_2,\mathcal{F}_2,P_2)\rightarrow(\widetilde{Y},\mathcal{B}(\widetilde{Y}))$
be two random elements as before. We assume $H_1[t]$ is independent
and identically distributed for different $t\in\mathbb{N}$ and so is
$H_2[t]$. Furthermore, $H_1$ is assumed to be independent with
$H_2$. The distribution functions of random variables $||H_1||$ and
$||H_2||$ are defined as $F_1(x):=P_1(||H_1||\le x)$ and
$F_2(y):=P_2(||H_2||\le y)$ for $x,y\in\mathbb{R}$, respectively. It
is worth noting that we do not ask the time scales $h_1$, $h_2$ to
be the same; and the coupled system may thus work in a kind of
asynchronous way.

Before going further, we give a definition for
\textit{nearly-emergence} that we adopt in this paper.
\smallskip

\noindent\textbf{Definition 1.}\itshape \quad Let $\mu,\nu>0$, $x\in
X$, $y\in\widetilde{Y}$, $\nu$(or $\mu$)-nearly-emergence occurs for
the population $\{1,\cdots,k\}$ if $||y||\le\nu$(or $||x||\le\mu$).
\normalfont\smallskip

Clearly, the exact emergence is no longer possible due to the random
perturbation.

\noindent\textbf{Remark 1.}\itshape \quad The two features $x$ and
$y$ of agents in the system I(D) are asymmetric and $y$ is the
object whose emergence behavior is of interest. In the system II(D),
the status of $x$ and $y$ is symmetric and both emergence behaviors
may be of interest. The same can be said for the continuous case in
Section 3 below.\normalfont

\medskip
\noindent{\large\textbf{2.2 Main results}}
\medskip

We define several constants that are only related with the initial
state $(x(0),y(0))$ of the population.

For system I(D):
$$
Q(\delta)=1\vee\frac1{\delta}, \quad
a=\frac{2C}GQ(\delta)||y(0)||^{\delta},\quad b=1+||x(0)||
$$
$$
\quad B_0=U_0-1,\quad
\mathcal{H}_0=\frac{2^{-\beta-1}G}{U_0^{\beta}},\quad
U_0=\left\{\begin{array}{ll}
\max\{(2a)^{\frac{1}{1-\gamma-\beta}},2b\},&
\mathrm{if}\ \beta+\gamma<1\\
\frac b{1-a},&\mathrm{if}\ \beta+\gamma=1\\
\frac{(\beta+\gamma)b}{\beta+\gamma-1},&\mathrm{if}\ \beta+\gamma>1
\end{array}
\right.
$$

For system II(D):
$$
\mathcal{H}_1=\frac{G_1}{2(1+||y(0)||)^{\beta_1}},\qquad
\mathcal{H}_2=\frac{G_2}{2(1+||x(0)||)^{\beta_2}}
$$

The main results in this section are stated as follows.

\smallskip
\noindent\textbf{Theorem 1.}\itshape \quad For dynamical system
I(D), we assume $ h<\min\Big\{\frac1{G},\
\frac1{2^{1-\gamma}C||y(0)||^{\delta}}\big(\frac{G}{2\mathcal{H}_0}\big)^{\frac{1-\gamma}{\beta}}\Big\},
$
and one of the following hypotheses holds:\\
(i) $\beta+\gamma<1$,\\
(ii)$\beta+\gamma=1$, and $||y(0)||<\big(\frac{G}{2CQ(\delta)}\big)^{\frac1{\delta}}$,\\
(iii)$\beta+\gamma>1$, and
$\big(\frac{1}{a(\beta+\gamma)}\big)^{\frac1{\beta+\gamma-1}}\frac{\beta+\gamma-1}{\beta+\gamma}>b+h\big((\frac{\beta+\gamma}{\beta+\gamma-1})b\big)^{\gamma}\frac{aG}{2Q(\delta)}$.\\
Then, for $\nu<||y(0)||$, $\nu$-nearly-emergence occurs in a number
of iterations bounded by $T_0:=\frac{2U_0^{\beta}}{hG}\ln
\big(\frac{||y(0)||}{\nu}\big)$ with probability at least
$F(\mathcal{H}_0\nu)^{T_0}$. In addition, if
$\mu<aU_0^{\beta+\gamma}$, let $T_1:=\frac{2U_0^{\beta}}{\delta
hG}\ln \big(\frac{aU_0^{\beta+\gamma}}{\mu}\big)$, then the events
$\{||x[t]-x[\tau]||\le\mu,\ for\ \tau>t\ge T_0\vee T_1\}$ and
$\{\nu$-nearly-emergence occurs in a number of iterations bounded by
$T_0\vee T_1\}$ hold simultaneously with probability at least
$F(\mathcal{H}_0\nu)^{T_0\vee T_1}$. \normalfont

\smallskip
\noindent\textbf{Theorem 2.}\itshape \quad For dynamical system
II(D), we assume $h_1<\frac1{G_1}$ and $h_2<\frac1{G_2}$. Then, for
$\mu<||x(0)||$, $\nu<||y(0)||$ with
$T_2:=\frac1{h_1\mathcal{H}_1}\ln \big(\frac{||x(0)||}{\mu}\big)$
and $T_3:=\frac1{h_2\mathcal{H}_2}\ln
\big(\frac{||y(0)||}{\nu}\big)$, $\mu$-nearly-emergence and
$\nu$-nearly-emergence both occur in a number of iterations bounded
by $T_2\vee T_3$ with probability at least
$(F_1(\mathcal{H}_1\mu)F_2(\mathcal{H}_2\nu))^{T_2\vee T_3}$.
\normalfont
\smallskip

We now give some concrete substances to illustrate emergence
behaviors of the general models I(D) and II(D).

For the system I(D), take $Y=\mathbb{R}^3$ with standard inner
product $\langle\cdot,\cdot\rangle$, and
$X=\widetilde{Y}=(\mathbb{R}^3)^k/\triangle\cong\triangle^{\perp}$.
For $u=(u_1,\cdots,u_k), v=(v_1,\cdots,v_k)\in\widetilde{Y}$, define
the inner product on $\widetilde{Y}$ as $\langle
u,v\rangle_{\widetilde{Y}}=\frac12\sum_{i,j=1}^{k}\langle
u_i-u_j,v_i-v_j \rangle$. Here $x\in X$ represents the spatial
positions of agents (e.g. birds, fishes, robots,...) and
$y\in\widetilde{Y}$ their velocities both projected to the subspace
$\triangle^{\perp}$ \cite{1}. Given $x\in X$, let the $k\times k$
matrix $A_x$ has entries $a_{ij}\ge\frac
K{(1+||x_i-x_j||)^{\beta}}$. Let $D_x$ be the $k\times k$ diagonal
matrix whose $i$th diagonal element is $d_i=\sum_{j\le k}a_{ij}$ and
$L_x:=D_x-A_x$. Then $L_x$ is the Laplacian of $A_x$ \cite{22}. For
$t\in\mathbb{N}$, take $J(x[t],y[t])=y[t]$, $S(x[t])=I_k-hL_x$, here
$I_k$ is the identity matrix of order $k$, and let the noise term
$H$ has the uniform distribution $U_{3k}(0,r)$ for some $r>0$ or the
Gaussian distribution $N(0,\sigma^2I_{3k})$ in the model I(D), and
then we recover the situations encountered in \cite{2}. Thm.1 in
\cite{2} is clearly a special case of Theorem 1 (and note that we
really said more). Other kinds of flocking scenarios such as
flocking with unrelated pairs and flocking with leader-follower
schemes can also be dealt with under our present framework (c.f.
\cite{3} Sect. 3). We mention here that the asymmetric conclusions
of $x$ and $y$ in Theorem 1 indeed give what we desire in a flocking
phenomenon; see \cite{1} (Rem. 2).

For the system II(D), let $\triangle_X$ be the diagonal of
$(\mathbb{R}^3)^k$ and take $X=(\mathbb{R}^3)^k/\triangle_X$ with
inner product defined as $\langle\cdot,\cdot\rangle_{\widetilde{Y}}$
above. Let $Y$ be the space of languages with some appropriate
distance defined on it (c.f. \cite{5}); and the metric of
$\widetilde{Y}$ is inherited from that of $Y$. Given $x\in
X,y\in\widetilde{Y}$, let $A_x=(a_{ij})$, $B_y=(b_{ij})$ be the
$k\times k$ matrices with entries $a_{ij}=f(||x_i-x_j||)$ and
$b_{ij}=g(||y_i-y_j||)$. $f,g:\mathbb{R}^+\rightarrow\mathbb{R}^+$
are some bounded non-increasing functions. For $t\in\mathbb{N}$,
take $S_1(y[t])=I_k-h_1L_1y$ and $S_2(x[t])=I_k-h_2L_2x$ with
$L_1y:=D_y-B_y$, $L_2x:=D_x-A_x$ in the model II(D). Here $D_x$ and
$D_y$ are $k\times k$ diagonal matrices defined similarly as above.
Computation of the distributions of $||H_1||$ and $||H_2||$ from
some proper random noises $H_1$, $H_2$ is a routine \cite{21}. Here
$x\in X$ is interpreted as the geographical positions of agents
projected to $\triangle_X^{\perp}$ and $y\in\widetilde{Y}$ as the
space of languages projected to $\triangle^{\perp}$. This
specification of system II(D) can be used to model emergence
behavior in linguistic evolution, since each agent tends to move to
others using similar languages and meanwhile the influence from
other agents' languages decreases according to distances \cite{1}
(Sect. 6).

\medskip
\noindent{\large\textbf{2.3 Proof of Theorem 1}}
\medskip

The proof closely follows that of \cite{2}, and we prove Theorem 1
through some intermediate steps.

\smallskip
\noindent\textbf{Proposition 1.}\itshape \quad Suppose
$0<h<\frac1G$, and $||H||\le\mathcal{H}_0||y[t]||$ for $0\le t<T$.
Then we have
$||y[t]||\le\big(1-\frac{hG}{(1+||x[t-1]||)^{\beta}}+h\mathcal{H}_0\big)||y[t-1]||$
for $t\le T$. Hence, for $1\le t\le T$,
$||y[t]||\le||y(0)||\prod_{i=0}^{t-1}\big(1-\frac{hG}{(1+||x[i]||)^{\beta}}+h\mathcal{H}_0\big)$.
\normalfont

\medskip
\noindent\textbf{Proof}. By (\ref{2}), we have
$$
||y[t]||=||S(x[t-1])y[t-1]+hH[t-1]||\le\Big(1-\frac{hG}{(1+||x[t-1]||)^{\beta}}+h\mathcal{H}_0\Big)||y[t-1]||.
$$ $\Box$

The following proposition is taken from \cite{4}.
\smallskip

\noindent\textbf{Proposition 2.}\itshape (\cite{4})\quad Let
$c_1,c_2>0$ and $s>q>0$. Then the equation $M(z)=z^s-c_1z^q-c_2=0$
has a unique positive zero $z^*$. Moreover,
$z^*\le\max\{(2c_1)^{\frac1{s-q}},(2c_2)^{\frac1s}\}$ and $M(z)\le0$
for $0\le z\le z^*$. \normalfont

\medskip
\noindent\textbf{Proposition 3.}\itshape \quad Let $T\in\mathbb{N}$
or $T=\infty$. Suppose $||H||\le\mathcal{H}_0||y[t]||$ for $0\le
t<T$ and $h<\min\Big\{\frac1{G},\
\frac1{2^{1-\gamma}C||y(0)||^{\delta}}\big(\frac{G}{2\mathcal{H}_0}\big)^{\frac{1-\gamma}{\beta}}\Big\}$,
and one of the following hypotheses holds:\\
(i) $\beta+\gamma<1$,\\
(ii)$\beta+\gamma=1$, and $||y(0)||<\big(\frac{G}{2CQ(\delta)}\big)^{\frac1{\delta}}$,\\
(iii)$\beta+\gamma>1$, and
$\big(\frac{1}{a(\beta+\gamma)}\big)^{\frac1{\beta+\gamma-1}}\frac{\beta+\gamma-1}{\beta+\gamma}>b+h\big((\frac{\beta+\gamma}{\beta+\gamma-1})b\big)^{\gamma}\frac{aG}{2Q(\delta)}$.\\
Then $1-\frac{hG}{2U_0^{\beta}}\in(0,1)$, for $0\le t<T$,
$||x[t]||\le B_0$ and $||y[t]||\le
||y(0)||\big(1-\frac{hG}{2U_0^{\beta}}\big)^t$. If $T=\infty$, then
$||y[t]||\rightarrow0$ as $t\rightarrow\infty$, and moreover, there
exists $\hat x\in X$ such that $x[t]\rightarrow\hat x$ as
$t\rightarrow\infty$, and $||x[t]-\hat x||\le
aU_0^{\beta+\gamma}\big(1-\frac{hG}{2U_0^{\beta}}\big)^{\delta t}$
for $t\ge0$. \normalfont

\medskip
\noindent\textbf{Proof}. Let $\Gamma:=\{0\le t\le T-1|\
(1+||x[t]||)^{\beta}\le\frac G{2\mathcal{H}_0}\}$. Assume
$\Gamma\ne\{0,\cdots,T-1\}$. Since $0\in\Gamma$ by the assumptions,
let $\hat t=\min\{\{0,\cdots,T-1\}\backslash\Gamma\}$. For $t<T$,
let $t^*$ be the point maximizing $||x||$ in $\{0,\cdots,t\}$. Then
by Proposition 1, for $t^*\le t<\hat t$ and $0\le i\le t$ we have
$$
\frac G{(1+||x[i]||)^{\beta}}-\mathcal{H}_0\ge\frac
G{(1+||x[t^*]||)^{\beta}}-\mathcal{H}_0\ge\frac
G{2(1+||x[^*]||)^{\beta}}:=R(t^*).
$$
By Proposition 1, for $\tau\le t$,
\begin{eqnarray*}
||x[\tau]||&\le&||x(0)||+\sum_{j=0}^{\tau-1}||x[j+1]-x[j]||\le||x(0)||+hC\sum_{j=0}^{\tau-1}(1+||x[j]||)^{\gamma}||y[j]||^{\delta}\\
 &\le&||x(0)||+||y(0)||^{\delta}hC\sum_{j=0}^{\tau-1}(1+||x[j]||)^{\gamma}\prod_{i=0}^{j-1}\big(1-\frac{hG}{(1+||x[i]||)^{\beta}}+h\mathcal{H}_0\big)^{\delta}\\
 &\le&||x(0)||+||y(0)||^{\delta}hC(1+||x[t^*]||)^{\gamma}\sum_{j=0}^{\tau-1}(1-hR(t^*))^{j\delta}\\
 &\le&||x(0)||+||y(0)||^{\delta}hC(1+||x[t^*]||)^{\gamma}\frac1{1-(1-hR(t^*))^{\delta}}\\
 &\le&||x(0)||+||y(0)||^{\delta}C(1+||x[t^*]||)^{\gamma}\frac{Q(\delta)}{R(t^*)}\\
 &=&||x(0)||+\frac{2CQ(\delta)}{G}||y(0)||^{\delta}(1+||x[t^*]||)^{\gamma+\beta}.
\end{eqnarray*}
Take $\tau=t^*$, and then we have
$$
1+||x[t^*]||\le1+||x(0)||+\frac{2CQ(\delta)}{G}||y(0)||^{\delta}(1+||x[t^*]||)^{\gamma+\beta}.
$$
Let $z=1+||x[t^*]||$, then the above inequality becomes
$M(z):=z-az^{\gamma+\beta}-b\le0$.

(i) Assume $\beta+\gamma<1$. By Proposition 2 and the definition of
$U_0$, we get $1+||x[t^*]||\le U_0$. Thereby for $t<\hat t$,
$||x[t]||\le U_0-1=B_0$. This implies
$$
(1+||x[t]||)^{\gamma+\beta}\le(1+||x[t^*]||)^{\gamma+\beta}\le
U_0^{\gamma+\beta}=\Big(\frac{2^{-\beta}G}{2\mathcal{H}_0}\Big)^{\frac{\gamma}{\beta}+1},
$$
for $t<\hat t$. Hence by Proposition 1 and $||y[t]||$ is decreasing,
\begin{eqnarray*}
1+||x[\hat
t]||&\le&1+||x[\hat{t}-1]||+hC(1+||x[\hat{t}-1]||)^{\gamma}||y[\hat{t}-1]||^{\delta}\\
 &\le&1+||x[\hat{t}-1]||+hC(1+||x[\hat{t}-1]||)^{\gamma}||y(0)||^{\delta}\\
 &\le&\Big(\frac{2^{-\beta}G}{2\mathcal{H}_0}\Big)^{\frac1{\beta}}+hC||y(0)||^{\delta}\Big(\frac{2^{-\beta}G}{2\mathcal{H}_0}\Big)^{\frac{\gamma}{\beta}}\le\Big(\frac{G}{2\mathcal{H}_0}\Big)^{\frac1{\beta}}.
\end{eqnarray*}
This is in contradiction with the definition of $\hat t$, therefore
$||x[t]||\le B_0$ and $\frac
G{(1+||x[t]||)^{\beta}}-\mathcal{H}_0\ge J_0:=\frac
G{2(1+B_0)^{\beta}}=\frac G{2U_0^{\beta}}$, for all $t<T$. It
follows from Proposition 1 that for $t<T$,
$$
||y[t]||\le||y(0)||\prod_{i=0}^{t-1}\big(1-\frac{hG}{(1+||x[i]||)^{\beta}}+h\mathcal{H}_0\big)\le||v(0)||(1-hJ_0)^t.
$$
For $T\ge\tau>t$, by employing (\ref{1}) we obtain
\begin{eqnarray*}
||x[\tau]-x[t]||&\le&\sum_{j=t}^{\tau-1}||x[j+1]-x[j]||=h\sum_{j=t}^{\tau-1}||J(x[j],y[j])||\\
 &\le&hC(1+B_0)^{\gamma}||v(0)||^{\delta}\sum_{j=t}^{\infty}(1-hJ_0)^{j\delta}=hCU_0^{\gamma}||v(0)||^{\delta}\frac{(1-hJ_0)^{t\delta}}{1-(1-hJ_0)^{\delta}}.
\end{eqnarray*}
When $T=\infty$, we take $t\rightarrow\infty$ to see that there
exists $\hat x\in X$ which is the limit of $x[t]$. Let
$\tau\rightarrow\infty$ in the above expression, we have
$$
||x[t]-\hat x||\le
CU_0^{\gamma}||v(0)||^{\delta}\frac{Q(\delta)}{J_0}(1-hJ_0)^{t\delta}\le
aU_0^{\beta+\gamma}\Big(1-\frac{hG}{2U_0^{\beta}}\Big)^{\delta t}.
$$

(ii) Assume $\beta+\gamma=1$. Then the inequality $M(z)\le0$ becomes
$$
(1+||x[t^*]||)\Big(1-\frac{2CQ(\delta)}{G}||y(0)||^{\delta}\Big)-(1+||x(0)||)\le0.
$$
By the assumption and definition of $B_0$, we get
$||x[t^*]||\le\frac{G(1+||x(0)||)}{G-2CQ(\delta)||y(0)||^{\delta}}-1=B_0$.
Thereby, for $t<\hat t$
$$
(1+||x[t]||)^{\gamma+\beta}\le1+||x[t^*]||\le\frac{G(1+||x(0)||)}{G-2CQ(\delta)||y(0)||^{\delta}}=\Big(\frac{2^{-\beta}G}{2\mathcal{H}_0}\Big)^{\frac1{\beta}}.
$$
We now proceed as in case (i).

(iii) Assume $\beta+\gamma>1$.
$M'(z)=1-(\beta+\gamma)az^{\beta+\gamma-1}$ has a unique zero
$z^*:=(\frac1{a(\gamma+\beta)})^{\frac1{\gamma+\beta-1}}$ and
$M(z^*)=(\frac1{a(\gamma+\beta)})^{\frac1{\gamma+\beta-1}}\frac{\gamma+\beta-1}{\gamma+\beta}-b>0$.
Since $M(0)=-b<0,M''(z)<0$ for $z>0$, $M$ is a convex function in
$(0,\infty)$ and has two positive zero $z_l$, $z_r$ satisfying
$0<z_l<z^*<z_r$.

For $t\in\mathbb{N}$, let $z(t):=1+||x[t^*]||$. We then have
$z(0)=1+||x(0)||=b<z^*$, which suggests $z(0)<z_l$. Now assume there
exists $t<\hat t$ such that $z(t)\ge z_r$ and denote $r$ be the
first such $t$. Hence, $r=r^*\ge1$ and $1+||x[t]||\le z(r-1)\le z_l$
for $t<r$. Let $z_0$ be the intersection of the $z$ axis with the
line joining $(0,-b)$ and $(z^*,M(z^*))$, (c.f. Fig.1 in \cite{2}).
Straightforward calculation yields $z_0=\frac
b{1-az^{*{\gamma+\beta-1}}}=\frac{(1+||x(0)||)(\gamma+\beta)}{\gamma+\beta-1}$.

By the definition of $B_0$, for $t<r$
\begin{equation}
||x[t]||\le z_l-1\le z_0-1=B_0\label{4}.
\end{equation}
We get $||x[r-1]||\le z_0-1$, $||x[r]||\ge z_r-1$ and accordingly,
$||x[r]-x[r-1]||\ge z^*-z_l$ since $z^*<z_r$. By the Lagrange
intermediate value theorem, there is $\xi\in[z_l,z^*]$ such that
$M(z^*)=M'(\xi)(z^*-z_l)$. Thereby $z^*-z_l\ge M(z^*)$ since $0\le
M'(\xi)\le 1$. Then $||x[r]-x[r-1]||\ge M(z^*)$. On the other hand,
$$
||x[r]-x[r-1]||\le hC(1+||x[r-1]||)^{\gamma}||y[r-1]||^{\delta}\le
hC(1+B_0)^{\gamma}||y(0)||^{\delta},
$$
which together with the above expression gives $M(z^*)\le
hCU_0^{\gamma}||y(0)||^{\delta}$. This is, however, contrary to our
hypothesis. Therefore, for all $t<\hat t$, $z(t)<z_l$ and
$||x[t]||\le B_0$ by noting inequality (\ref{4}). We then obtain,
$$
(1+||x[t]||)^{\gamma+\beta}\le z_0^{\gamma+\beta}=
\Big(\frac{(1+||x(0)||)(\gamma+\beta)}{\gamma+\beta-1}\Big)^{\gamma+\beta}=\Big(\frac{2^{-\beta}G}{2\mathcal{H}_0}\Big)^{\frac{\gamma}{\beta}+1}
$$
for $t<\hat t$. Now we proceed as in case (i). $\Box$

\medskip
\noindent\textbf{Proof of Theorem 1}. Suppose the conditions of
Proposition 3 hold for some $T>0$. Then we have $||y[t]||\le
||y(0)||\big(1-\frac{hG}{2U_0^{\beta}}\big)^t$ for $t<T$ and
$||x[\tau]-x[t]||\le
aU_0^{\beta+\gamma}\big(1-\frac{hG}{2U_0^{\beta}}\big)^{\delta t}$
for $t<\tau\le T$. If $T=\infty$, $||x[t]-\hat x||\le
aU_0^{\beta+\gamma}\big(1-\frac{hG}{2U_0^{\beta}}\big)^{\delta t}$
for $t\ge0$. By the proof of Proposition 3 and straightforward
calculations, we have $||y[T]||\le\nu$ when $T\ge T_0$;
$||x[\tau]-x[t]||\le\mu$ when $\tau>t\ge T_1$; and $||x[t]-\hat
x||\le\mu$ when $t\ge T_1$. If $\nu$-nearly-emergence has not
occurred, then $||y[t]||\ge\nu$, wherefore by the definition of
function $F$,
$$
P(||H[t]||\le\mathcal{H}_0||y[t]||)\ge
P(||H[t]||\le\mathcal{H}_0\nu)=F(\mathcal{H}_0\nu).
$$
Since $\{H[t]\}$ are i.i.d. for varying $t$,  we get
$$
P(||H[t]||\le\mathcal{H}_0||y[t]||\ \mathrm{for}\
t=0,\cdots,T_0-1)\ge F(\mathcal{H}_0\nu)^{T_0},
$$
which yields the first part of the conclusions. Likewise, we have
$$
P(||H[t]||\le\mathcal{H}_0||y[t]||\ \mathrm{for}\ t=0,\cdots,T_0\vee
T_1-1)\ge F(\mathcal{H}_0\nu)^{T_0\vee T_1}.
$$
We then conclude the proof. $\Box$

\medskip
\noindent{\large\textbf{2.4 Proof of Theorem 2}}
\medskip

\smallskip
\noindent\textbf{Proposition 4.}\itshape \quad Suppose
$0<h_1<\frac1{G_1}$, $0<h_2<\frac1{G_2}$, and
$||H_1[t_1]||\le\mathcal{H}_1||x[t_1]||$ for $0\le t_1<T^1$;
$||H_2[t_2]||\le\mathcal{H}_2||y[t_2]||$ for $0\le t_2<T^2$. Then we
have
$||x[t_1]||\le\big(1-\frac{h_1G_1}{(1+||y[t_1-1]||)^{\beta_1}}+h_1\mathcal{H}_1\big)||x[t_1-1]||$
for $t_1\le T^1$ and
$||y[t_2]||\le\big(1-\frac{h_2G_2}{(1+||x[t_2-1]||)^{\beta_2}}+h_2\mathcal{H}_2\big)||y[t_2-1]||$
for $t_2\le T^2$. Consequently, for $t_1\le T^1$,
$||x[t_1]||\le||x(0)||\prod_{i=0}^{t-1}\big(1-\frac{h_1G_1}{(1+||y[i]||)^{\beta_1}}+h_1\mathcal{H}_1\big)$;
and for $t_2\le T^2$,
$||y[t_2]||\le||y(0)||\prod_{j=0}^{t-1}\big(1-\frac{h_2G_2}{(1+||x[j]||)^{\beta_2}}+h_2\mathcal{H}_2\big)$.
\normalfont

\medskip
\noindent\textbf{Proof}. The proof readily follows by utilizing
condition (\ref{3}). $\Box$
\smallskip

The following proposition is critical to the proof of Theorem 2.

\smallskip
\noindent\textbf{Proposition 5.}\itshape \quad Let
$T^1,T^2\in\mathbb{N}\cup\{\infty\}$. Suppose
$||H_1||\le\mathcal{H}_1||x[t_1]||$ for $0\le t_1<T^1$;
$||H_2||\le\mathcal{H}_2||y[t_2]||$ for $0\le t_2<T^2$ and
$h_1<\frac1{G_1}$, $h_2<\frac1{G_2}$. Then, for $0\le t<T^1\wedge
T^2$, $||x[t]||\le||x(0)||(1-h_1\mathcal{H}_1)^{t}$ and
$||y[t]||\le||y(0)||(1-h_2\mathcal{H}_2)^{t}$. \normalfont

\medskip
\noindent\textbf{Proof}. Denote $T:=T^1\wedge T^2$ for brevity. Let
$\Gamma_1:=\{0\le t_1\le T-1|\
(1+||y[t_1]||)^{\beta_1}\le\frac{G_1}{2\mathcal{H}_1}\}$. Assume
$\Gamma_1\ne\{0,\cdots,T-1\}$. Since $0\in\Gamma_1$ by the
assumptions, let $\hat
t_1=\min\{\{0,\cdots,T-1\}\backslash\Gamma_1\}$. For $t_1<T$, let
$t_1^*$ be the point maximizing $||y||$ in $\{0,\cdots,t_1\}$.
Analogously, we define $\Gamma_2:=\{0\le t_2\le T-1|\
(1+||x[t_2]||)^{\beta_2}\le\frac{G_2}{2\mathcal{H}_2}\}$ and assume
$\Gamma_2\ne\{0,\cdots,T-1\}$. Let $\hat
t_2=\min\{\{0,\cdots,T-1\}\backslash\Gamma_2\}$. For $t_2<T$, let
$t_2^*$ be the point maximizing $||x||$ in $\{0,\cdots,t_2\}$.

Then by Proposition 4, for $t_1^*\le t_1<\hat t_1$ and $0\le i\le
t_1$ we have
$$
\frac{G_1}{(1+||y[i]||)^{\beta_1}}-\mathcal{H}_1\ge\frac
{G_1}{(1+||y[t_1^*]||)^{\beta_1}}-\mathcal{H}_1\ge\frac
{G_1}{2(1+||y[t_1^*]||)^{\beta_1}}:=R_1(t_1^*).
$$
By exploiting Proposition 4, for $\tau\le t_1$,
\begin{eqnarray}
||x[\tau]||&\le&||x(0)||\prod_{i=0}^{\tau-1}\big(1-\frac{h_1G_1}{(1+||y[i]||)^{\beta_1}}+h_1\mathcal{H}_1\big)\nonumber\\
 &\le&||x(0)||\prod_{i=0}^{\tau-1}(1-h_1R_1(t_1^*))=||x(0)||(1-h_1R_1(t_1^*))^{\tau}.\label{5}
\end{eqnarray}
Likewise, for $t_2^*\le t_2<\hat t_2$ and $0\le j\le t_2$ we have
$$
\frac{G_2}{(1+||x[j]||)^{\beta_2}}-\mathcal{H}_2\ge\frac
{G_2}{(1+||x[t_2^*]||)^{\beta_2}}-\mathcal{H}_2\ge\frac
{G_2}{2(1+||x[t_2^*]||)^{\beta_2}}:=R_2(t_2^*).
$$
And for $\tau\le t_2$, we have
\begin{equation}
||y[\tau]||\le||y(0)||\prod_{j=0}^{\tau-1}\big(1-\frac{h_2G_2}{(1+||x[j]||)^{\beta_2}}+h_2\mathcal{H}_2\big)
\le||y(0)||(1-h_2R_2(t_2^*))^{\tau}.\label{6}
\end{equation}

Now we shall compare $\hat t_1$ with $\hat t_2$ to deduce
contradictions. If $\hat t_1<\hat t_2$. We take $t_2=\hat t_1$ and
get from (\ref{6}), $||y[\hat t_1]||\le
||y(0)||(1-h_2R_2(t_2^*))^{\hat t_1}\le||y(0)||$. This is in
contradiction with definition of $\hat t_1$. Similarly, if $\hat
t_1>\hat t_2$, we can also deduce a contradiction. Now let's
consider the case $\hat t_1=\hat t_2$. If so, take $\tau=t_2^*$ in
(\ref{5}) and $\tau=t_1^*$ in (\ref{6}). We then have
$||x[t_2^*]||\le||x(0)||(1-h_1R_1(t_2^*))^{t_2^*}$ and
$||y[t_1^*]||\le||y(0)||(1-h_2R_2(t_1^*))^{t_1^*}$. These
inequalities imply $t_1^*=t_2^*=0$. That is to say, for all $t<\hat
t_1=\hat t_2$, $||x[t]||\le||x(0)||$ and $||y[t]||\le||y(0)||$.
Combining these with Proposition 4, we derive
\begin{eqnarray*}
||y[\hat t_1]||&\le&\big(1-\frac{h_2G_2}{(1+||x[\hat
t_1-1]||)^{\beta_2}}+h_2\mathcal{H}_2\big)||y[\hat t_1-1]||\\
 &=&\big(1-\frac{h_2G_2}{(1+||x[\hat
t_1-1]||)^{\beta_2}}+\frac{h_2G_2}{(1+||x(0)||)^{\beta_2}}\big)||y[\hat t_1-1]||\\
 &\le&||y[\hat t_1-1]||\le||y(0)||.
\end{eqnarray*}
Here is the contradiction as in the first case above. It is at this
stage that we may conclude that $\hat t_1$ and $\hat t_2$ must not
\textit{both} exist.

Without loss of generality, we assume $\hat t_1$ does not exist,
i.e. $\Gamma_1=\{0,\cdots,T-1\}$. For $t<T$,
$(1+||y[t]||)^{\beta_1}\le\frac{G_1}{2\mathcal{H}_1}=(1+||y(0)||)^{\beta_1}$.
Hence, $||y[t]||\le||y(0)||$ and
$R_1(t_1^*)\ge\frac{G_1}{2(1+||y(0)||)^{\beta_1}}$. From (\ref{5}),
we have $||x[t]||\le||x(0)||(1-h_1\mathcal{H}_1)^{t}$. If $\hat t_2$
exists, we then take $t=\hat t_2$ in the above inequality to get
$||x[\hat t_2]||\le||x(0)||$. It is in contradiction with the
definition of $\hat t_2$. Therefore, $\hat t_2$ does not exist
either. Now we can proceed as above to get
$||y[t]||\le||y(0)||(1-h_2\mathcal{H}_2)^{t}$ for $t<T$, by using
(\ref{6}). $\Box$

\medskip
\noindent\textbf{Proof of Theorem 2}. Suppose the conditions of
Proposition 5 hold for some $T^1,T^2>0$. Then we have
$||x[t]||\le||x(0)||(1-h_1\mathcal{H}_1)^{t}$ and
$||y[t]||\le||y(0)||(1-h_2\mathcal{H}_2)^{t}$ for $t<T:=T^1\wedge
T^2$. By direct calculations, we get $||x[T]||\le\mu$ when $T\ge
T_2$ and $||y[T]||\le\nu$ when $T\ge T_3$. If $\mu$-nearly-emergence
and $\nu$-nearly-emergence have not occurred at time $t_1h_1$ and
$t_2h_2$ resp., then $||x[t_1]||\ge\mu$ and $||y[t_2]||\ge\nu$.
Whence,
$$
P_1(||H_1[t_1]||\le\mathcal{H}_1||x[t_1]||)\ge
P_1(||H_1[t_1]||\le\mathcal{H}_1\mu)=F_1(\mathcal{H}_1\mu),
$$
and
$$
P_2(||H_2[t_2]||\le\mathcal{H}_2||y[t_2]||)\ge
P_2(||H_2[t_2]||\le\mathcal{H}_2\nu)=F_2(\mathcal{H}_2\nu).
$$
Since $\{H_i[t_i]\}$ are i.i.d. for varying $t_i$, $i=1,2$, and
$H_1$ is independent with $H_2$, by letting $P=P_1\times P_2$ be the
independent product of $P_1$ and $P_2$ (c.f. \cite{21}), we get
\begin{multline}
P(||H_1[t]||\le\mathcal{H}_1||x[t]||\ \mathrm{and}\
||H_2[t]||\le\mathcal{H}_2||y[t]||,\ \mathrm{for}\
t=0,\cdots,T_2\vee
T_3-1)\\
\ge (F_1(\mathcal{H}_1\mu)F_2(\mathcal{H}_2\nu))^{T_2\vee
T_3},\nonumber
\end{multline}
which concludes the proof. $\Box$

\medskip
\noindent{\Large\textbf{3. Continuous-time emergence}}
\smallskip

\noindent{\large\textbf{3.1 Models setup (I(C), II(C))}}
\medskip

In principle, by letting the time steps $h$, $h_i$ approach zero, we
may derive the continuous counterparts of systems I(D) and II(D). We
shall, however, make some modifications for technical reason and it
is at this time the inner product structures of spaces
$X,Y,\widetilde{Y}$ take effect.

Let $t\in\mathbb{R}^+$, for $x\in X$, $y\in \widetilde{Y}$, consider
the following dynamical system:
$$
I(C):\qquad\Big\{\begin{array}{l} x'(t)=J(x(t),y(t))\\
y'(t)=-L_xy(t)+hH(t)
\end{array}
$$
Here, as in the system I(D), $J:X\times\widetilde{Y}\rightarrow X$
is a Lipschitz or $C^1$ operator. We now require, for some
$C,\delta>0$, $0\le \gamma<1$, that
\begin{equation}
||J(x,y)||\le C(1+||x||^2)^{\frac{\gamma}{2}}||y||^{\delta}\label{7}
\end{equation}
for all $x\in X$, $y\in\widetilde{Y}$. Denote $\mathbb{R}^{k\times
k}$ as the space of $k\times k$ real matrices. Let $L:
X\rightarrow\mathbb{R}^{k\times k}$ be a Lipschitz or $C^1$ operator
with $L: x\mapsto L_x$. $L_x$ can be seen as a linear transformation
on $Y^k$ by mapping $(y_1,\cdots,y_k)\in Y^k$ to
$(L_x(i,1)y_1+\cdots+L_x(i,k)y_k)_{i\le k}$. Here, $L_x(i,j)$ is the
$(i,j)$ entry of matrix $L_x$. For $x\in X$, define
$\phi_x:=\min_{y\ne0 \atop y\in\widetilde{Y}}\frac{\langle
L_xy,y\rangle}{||y||^2}$. We impose the following two hypotheses on
$L$: (i) For $x\in X$, $y\in Y$, $L_x(y,\cdots,y)=0$. (ii) there
exists $K>0,\beta\ge0$ such that
\begin{equation}
\phi_x\ge\frac{K}{(1+||x||^2)^{\beta}}\label{8}
\end{equation}
for all $x\in X$. It is easy to see from (i) that $L_x$ induces a
linear transformation on $\widetilde{Y}$, which will also be denoted
as $L_x$ for notational simplicity. Let $H(t)$ be a continuous time
stochastic process defined on some probability space
$(\Omega,\mathcal{F},P)$ taking value in
$(\widetilde{Y},\mathcal{B}(\widetilde{Y}))$. Let
$F:\mathbb{R}^+\times\mathbb{R}^+\rightarrow[0,1]$ be a real
function (not necessarily a distribution function) such that
$P(\max_{0\le t\le T}||H(t)||\le x)\ge F(x,T)$, for
$x,T\in\mathbb{R}^+$. We may observe that $F(x,T)$ is non-decreasing
w.r.t. $x$ while non-increasing w.r.t. $T$.

Next, we introduce a continuous version of II(D). For $x\in X$,
$y\in \widetilde{Y}$, consider the following dynamical system:
$$
II(C):\qquad\Big\{\begin{array}{l} x'(t)=-L_{1y}x(t)+H_1(t)\\
y'(t)=-L_{2x}y(t)+H_2(t)
\end{array}
$$
Here operators $L_1,L_2$ are given similarly as $L$ above with $L_1$
defined on $\widetilde{Y}$ and $L_2$ on $X$. For $x\in X$, $y\in
\widetilde{Y}$, define $\xi_x:=\min_{y\ne0 \atop
y\in\widetilde{Y}}\frac{\langle L_{2x}y,y\rangle}{||y||^2}$ and
$\eta_y:=\min_{x\ne0 \atop x\in X}\frac{\langle
L_{1y}x,x\rangle}{||x||^2}$. The corresponding hypothesis (ii) above
becomes (ii$'$): there exist $K_1,K_2>0$, $\beta_1,\beta_2\ge0$ such
that
\begin{equation}
\xi_x\ge\frac{K_1}{(1+||x||^2)^{\beta_1}},\qquad
\eta_y\ge\frac{K_2}{(1+||y||^2)^{\beta_2}}\label{9}
\end{equation}
for all $x\in X$, $y\in\widetilde{Y}$. Let
$H_1(t):(\Omega_1,\mathcal{F}_1,P_1)\rightarrow(X,\mathcal{B}(X))$
and
$H_2(t):(\Omega_2,\mathcal{F}_2,P_2)\rightarrow(\widetilde{Y},\mathcal{B}(\widetilde{Y}))$
be two continuous time stochastic processes, which are independent
with each other. Let $i=1,2$,
$F_i:\mathbb{R}^+\times\mathbb{R}^+\rightarrow[0,1]$ be real
functions such that $P_i(\max_{0\le t\le T}||H_i(t)||\le x)\ge
F_i(x,T)$, for $x,T\in\mathbb{R}^+$.

For notational convenience we sometimes write $L_t:=L_{x(t)}$,
$\phi_t:=\phi_{x(t)}$, $\xi_t:=\xi_{x(t)}$ and
$\eta_t:=\eta_{y(t)}$.

\medskip
\noindent{\large\textbf{3.2 Main results}}
\medskip

As in the discrete case, we define several constants which are only
dependent on the initial state $(x(0),y(0))$ of the population.

For system I(C):
$$
a=2^{\frac{1+\gamma+2\beta}{1-\gamma}}\frac{((1-\gamma)C)^{\frac2{1-\gamma}}||y(0)||^{\frac{2\delta}{1-\gamma}}}{(\delta
K)^{\frac2{1-\gamma}}},\quad
b=2^{\frac{1+\gamma}{1-\gamma}}(1+||x(0)||^2),\quad
\alpha=\frac{2\beta}{1-\gamma}
$$
$$
\quad B_0=U_0-1,\quad
B_1=\frac{2C||y(0)||^{\delta}B_0^{\frac{\gamma}2+\beta}}{\delta
K},\quad \mathcal{H}_0=\frac{2^{-\beta-1}G}{U_0^{\beta}}
$$
$$
U_0=\left\{\begin{array}{ll}
\max\{(2a)^{\frac{1-\gamma}{1-\gamma-2\beta}},2b\},&
\mathrm{if}\ 2\beta+\gamma<1\\
\frac b{1-a},&\mathrm{if}\ 2\beta+\gamma=1\\
(\frac1{a\alpha})^{\frac1{\alpha-1}},&\mathrm{if}\ 2\beta+\gamma>1
\end{array}
\right.
$$

For system II(C):
$$
\mathcal{H}_1=\frac{K_2}{2(1+||y(0)||^2)^{\beta_2}},\qquad
\mathcal{H}_2=\frac{K_1}{2(1+||x(0)||^2)^{\beta_1}}
$$

The main results in this section are stated as follows.

\smallskip
\noindent\textbf{Theorem 3.}\itshape \quad Let $x(0)\in X$ and
$y(0)\in \widetilde{Y}$, then there exists a unique solution
$(x(t),y(t))$ of the dynamical system I(C) for all $t\in\mathbb{R}$.
Moreover, assume one of the following hypotheses holds:\\
(i) $2\beta+\gamma<1$,\\
(ii)$2\beta+\gamma=1$, and $||y(0)||<\big(\frac{(\delta K)^2}{2^{1+\gamma+2\beta}((1-\gamma)C)^2}\big)^{\frac1{2\delta}}$,\\
(iii)$2\beta+\gamma>1$, and
$\big(\frac{1}{a\alpha}\big)^{\frac1{\alpha-1}}\frac{\alpha-1}{\alpha}>b$.\\
Then, for $\nu<||y(0)||$, $\nu$-nearly-emergence occurs before time
$T_0:=\frac{2B_0^{\beta}}{K}\ln \big(\frac{||y(0)||}{\nu}\big)$ with
probability at least $F(\mathcal{H}_0\nu,T_0)$. In addition, if
$\mu<B_1$, let $T_1:=\frac{2B_0^{\beta}}{K\delta}\ln
\big(\frac{B_1}{\mu}\big)$, then the events
$\{||x[t]-x[\tau]||\le\mu,\ for\ \tau>t\ge T_0\vee T_1\}$ and
$\{\nu$-nearly-emergence occurs before time $T_0\vee T_1\}$ hold
simultaneously with probability at least $F(\mathcal{H}_0\nu,T_0\vee
T_1)$. \normalfont

\smallskip
As in Section 2.2, we may readily recover the continuous-time result
in \cite{2} by letting $L_x$ be the Laplacian of $A_x$, for $x\in
X$; and taking $J(x,y)=y$ and the coordinate processes of $H(t)$ as
independent smoothed Wiener processes.

\smallskip
\noindent\textbf{Theorem 4.}\itshape \quad Let $x(0)\in X$ and
$y(0)\in \widetilde{Y}$, then there exists a unique solution
$(x(t),y(t))$ of the dynamical system II(C) for all
$t\in\mathbb{R}$. Furthermore, for $\mu<||x(0)||$, $\nu<||y(0)||$
with $T_2:=\frac{2(1+||y(0)||^2)^{\beta_2}}{K_2}\ln
\big(\frac{||x(0)||}{\mu}\big)$ and
$T_3:=\frac{2(1+||x(0)||^2)^{\beta_1}}{K_1}\ln
\big(\frac{||y(0)||}{\nu}\big)$, either $\mu$-nearly-emergence or
$\nu$-nearly-emergence occurs before time $T_2\vee T_3$ with
probability at least\\ $F_1(\mathcal{H}_1\mu,T_2\vee
T_3)F_2(\mathcal{H}_2\nu,T_2\vee T_3)$. \normalfont

\smallskip
Compared with Theorem 2, the last result is weaker due to the fact
that in the continuous case the stochastic processes $H_i(t)$ do not
possess ``independence property'' among different ``time steps''.
However, if the noise does not impose on both equations of system
II(C), we have the following corollary.

\smallskip
\noindent\textbf{Corollary 1.}\itshape\quad Suppose $H_1(t)\equiv0$.
Under the assumptions of Theorem 4, the events
$\{\mu$-nearly-emergence occurs before time $T_2\}$ and
$\{\nu$-nearly-emergence occurs before time $T_3\}$ hold
simultaneously with probability at least
$F_2(\mathcal{H}_2\nu,T_3)$. An analogous result holds for the case
$H_2(t)\equiv0$. \normalfont

We mention that it is possible to have results similar with
Corollary 1 when $H_1(t)$ is small enough or possesses independent
increments.

\medskip
\noindent{\large\textbf{3.3 Proof of Theorem 3}}
\medskip

The main procedure of the proof follows that of \cite{2}. Let
$\theta_t:=\min_{\tau\in[0,t]}\phi_{\tau}$,
$\Gamma(t):=\Gamma(x(t)):=||x(t)||^2$ and
$\Lambda(t):=\Lambda(y(t)):=||y(t)||^2$.

\smallskip
\noindent\textbf{Proposition 6.}\itshape \quad Suppose for $0\le
t<T$, $||H(t)||\le ||y(t)||\mathcal{H}_0$. Then we have
$\Lambda(t)\le \Lambda(0)e^{-2t(\theta_t-\mathcal{H}_0)}$ for $0\le
t<T$. \normalfont

\medskip
\noindent\textbf{Proof}. The proof parallels with that of \cite{2}
(Prop. 4), by using assumption (\ref{8}). Hence we omit it. $\Box$

\smallskip
\noindent\textbf{Proposition 7.}\itshape \quad Suppose for $0\le
t<T$, $\theta_t>\mathcal{H}_0$. Then we have
$$
\Gamma(t)\le2^{\frac{1+\gamma}{1-\gamma}}\Big((1+\Gamma(0))+
\frac{((1-\gamma)C)^{\frac2{1-\gamma}}\Lambda(0)^{\frac{\delta}{1-\gamma}}}{(\delta(\theta_t-\mathcal{H}_0))^{\frac2{1-\gamma}}}
\Big)-1,
$$
for $0\le t<T$. \normalfont

\medskip
\noindent\textbf{Proof}. Mimicking the proof of \cite{3} (Prop. 6)
by employing Proposition 6. We leave this proof as an exercise for
the reader. $\Box$

\smallskip
\noindent\textbf{Proposition 8.}\itshape \quad Suppose for $0\le
t<T$, $||H(t)||\le ||y(t)||\mathcal{H}_0$. Assume one of the
following hypotheses holds:\\
(i) $2\beta+\gamma<1$,\\
(ii)$2\beta+\gamma=1$, and $\Lambda(0)<\big(\frac{(\delta K)^2}{2^{1+\gamma+2\beta}((1-\gamma)C)^2}\big)^{\frac1{\delta}}$,\\
(iii)$2\beta+\gamma>1$, and
$\big(\frac{1}{a\alpha}\big)^{\frac1{\alpha-1}}\frac{\alpha-1}{\alpha}>b$.\\
Then for $0\le t<T$, $\Gamma(t)\le B_0$ and $\Lambda(t)\le
\Lambda(0)e^{-\frac{Kt}{B_0^{\beta}}}$. If $T=\infty$, then
$\Lambda(t)\rightarrow0$ as $t\rightarrow\infty$ and there exists
$\hat x\in X$ such that $x(t)\rightarrow\hat x$ as
$t\rightarrow\infty$ and $||x(t)-\hat x||\le B_1e^{-\frac{\delta
Kt}{2B_0^{\beta}}}$ for $t\ge0$. \normalfont

\medskip
\noindent\textbf{Proof}. Let $\Xi:=\{t\in[0,T)|\
(1+||x(t)||^2)^{\beta}\le\frac K{2\mathcal{H}_0}\}$. Assume
$\Xi\ne[0,T)$. Since $0\in\Xi$ by the assumptions, let $\hat
t=\inf\{[0,T)\backslash\ \Xi\}$. Let $t<\hat t$ and $t^*\in[0,t]$ be
the point maximizing $||x||$ in $[0,t]$. Then by Proposition 6,
$$
\theta_t=\min_{\tau\in[0,t]}\phi_{\tau}\ge\min_{\tau\in[0,t]}\frac
K{(1+\Gamma(\tau))^{\beta}}=\frac K{(1+\Gamma(t^*))^{\beta}}.
$$
Since $t^*\le t<\hat t$, $t^*\in\Xi$, we deduce
$$
\theta_t-\mathcal{H}_0\ge\frac
K{(1+\Gamma(t^*))^{\beta}}-\mathcal{H}_0\ge\frac
K{2(1+\Gamma(t^*))^{\beta}}>0.
$$
Thereby Proposition 7 implies,
$$
\Gamma(t)\le2^{\frac{1+\gamma}{1-\gamma}}\Big((1+\Gamma(0))+
\frac{((1-\gamma)C)^{\frac2{1-\gamma}}\Lambda(0)^{\frac{\delta}{1-\gamma}}}{(\delta
K)^{\frac2{1-\gamma}}}\big(2(1+\Gamma(t^*))\big)^{\frac{2\beta}{1-\gamma}}\Big)-1.
$$
Take $t=t^*$ in the above expression, and then we have
\begin{equation}
(1+\Gamma(t^*))-2^{\frac{1+\gamma}{1-\gamma}}
\frac{((1-\gamma)C)^{\frac2{1-\gamma}}\Lambda(0)^{\frac{\delta}{1-\gamma}}}{(\delta
K)^{\frac2{1-\gamma}}}\big(2(1+\Gamma(t^*))\big)^{\frac{2\beta}{1-\gamma}}-2^{\frac{1+\gamma}{1-\gamma}}(1+\Gamma(0))\le0.\label{10}
\end{equation}
Let $z=1+\Gamma(t^*)$. Then inequality (\ref{10}) becomes $M(z)\le0$
with $M(z):=z-az^{\alpha}-b$. To complete the proof, we may argue as
in Proposition 3 dividing into three cases and consult the proof of
Thm. 2 of \cite{3}. $\Box$

\medskip
\noindent\textbf{Proof of Theorem 3}. The existence of a unique
solution for each $\omega\in\Omega$ follows from \cite{8} (Ch. 8).

By assumption, we have for $\varepsilon>0$, $P(\max_{0\le t\le
T}||H(t)||\le \varepsilon)\ge F(\varepsilon,T)$. Define
$T(\omega):=\inf\{t\ge0|\ ||y(t)||\le\nu\}$ and take
$\varepsilon=\mathcal{H}_0\nu$. To prove the first claim in Theorem
3, it suffices to proof
\begin{equation}
P(T(\omega)\le T_0)\ge F(\mathcal{H}_0\nu,T_0).\label{11}
\end{equation}

Take $\omega\in\{\omega\in\Omega|\ \max_{0\le t\le
T_0}||H(t)||\le\mathcal{H}_0\nu\}$. Assume $T(\omega)>T_0$. Then
there is some $\zeta>0$ such that $T(\omega)>T_0+\zeta$. Wherefore,
$||y(t)||>\nu$ on the interval $[0,T_0+\zeta]$. In particular,
$||y(T_0)||>\nu$. Since $\max_{0\le t\le
T_0}||H(t)||\le\mathcal{H}_0\nu<||y(t)||\mathcal{H}_0$, Proposition
8 holds for $T=T_0$. Then for $0\le t<T_0$, we have $||y(t)||\le
||y(0)||e^{-\frac{Kt}{2B_0^{\beta}}}$. As $y(t)$ is continuous, we
obtain $||y(T_0)||\le ||y(0)||e^{-\frac{KT_0}{2B_0^{\beta}}}=\nu$.
This is a contradiction, which then finishes the proof of
(\ref{11}).

Next, by the proof of Proposition 8 and straightforward calculations
(much the same way as in the proof of Theorem 1), we have
$||x(\tau)-x(t)||\le\mu$ when $\tau>t\ge T_1$; and $||x(t)-\hat
x||\le\mu$ when $t\ge T_1$. We then conclude the proof by replace
$T_0$ with $T_0\vee T_1$ in the above argument. $\Box$

\medskip
\noindent{\large\textbf{3.4 Proof of Theorem 4}}
\medskip

Denote $\Lambda_1(x(t)):=\langle x(t),x(t) \rangle$ and
$\Lambda_2(y(t)):=\langle y(t),y(t) \rangle$. In the sequel, we will
suppress the subscript $1,2$ of $\Lambda$ for brevity.

\smallskip
\noindent\textbf{Proposition 9.}\itshape \quad Suppose
$||H_1(t_1)||\le||x(t_1)||\mathcal{H}_1$ for $0\le t_1<T^1$ and
$||H_2(t_2)||\le||y(t_2)||\mathcal{H}_2$ for $0\le t_2<T^2$. Then
for $0\le t_1<T^1$,
$\Lambda(x(t_1))\le\Lambda(x(0))e^{-2\int_0^{t_1}\eta_{\tau}-\mathcal{H}_1\mathrm{d}\tau}$;
and for $0\le t_2<T^2$,
$\Lambda(y(t_2))\le\Lambda(y(0))e^{-2\int_0^{t_2}\xi_{\tau}-\mathcal{H}_2\mathrm{d}\tau}$.
\normalfont

\medskip
\noindent\textbf{Proof}. Let $\tau\in[0,t_1]$. Therefore, by using
(\ref{9}) we have
\begin{eqnarray*}
\Lambda'(x(\tau))&=&\langle x(\tau),x(\tau)\rangle=2\langle
x'(\tau),x(\tau)\rangle=-2\langle L_yx(\tau)-H_1(\tau),x(\tau)\rangle\\
 &\le&-2\eta_{\tau}\Lambda(x(\tau))+2||H_1(\tau)||\cdot||x(\tau)||=-2\Lambda(x(\tau))(\eta_{\tau}-\mathcal{H}_1)\\
\end{eqnarray*}
Integrating the above inequality from $0$ to $t_1$ deduces the first
result. The other can be proved likewise. $\Box$

\smallskip
\noindent\textbf{Proposition 10.}\itshape \quad Suppose
$||H_1(t_1)||\le\mathcal{H}_1||x(t_1)||$ for $0\le t_1<T^1$; and
$||H_2(t_2)||\le\mathcal{H}_2||y(t_2)||$ for $0\le t_2<T^2$. Then we
obtain for $0\le t<T^1\wedge T^2$,
$\Lambda(x(t))\le\Lambda(x(0))e^{-2t\mathcal{H}_1}$ and
$\Lambda(y(t))\le\Lambda(y(0))e^{-2t\mathcal{H}_2}$. \normalfont

\medskip
\noindent\textbf{Proof}. Denote $T:=T^1\wedge T^2$. Let
$\Gamma_1:=\{t_1\in[0,T)|\
(1+\Lambda(y(t_1)))^{\beta_2}\le\frac{K_2}{2\mathcal{H}_1}\}$.
Assume $\Gamma_1\ne[0,T)$. Since $0\in\Gamma_1$ by the assumptions,
let $\hat t_1=\inf\{[0,T)\backslash\Gamma_1\}$. Let $t_1<\hat t_1$
and $t_1^*\in[0,t_1]$ be the point maximizing $\Lambda(y(t))$ in
$[0,t_1]$. Analogously, we define $\Gamma_2:=\{t_2\in[0,T)|\
(1+\Lambda(x(t_2)))^{\beta_1}\le\frac{K_1}{2\mathcal{H}_2}\}$ and
assume $\Gamma_2\ne[0,T)$. Let $\hat
t_2=\inf\{[0,T)\backslash\Gamma_2\}$. Let $t_2<\hat t_2$ and
$t_2^*\in[0,t_2]$ be the point maximizing $\Lambda(x(t))$ in
$[0,t_2]$.

Since $t_1^*\le t_1<\hat t_1$, $t_1^*\in\Gamma_1$, and then we have
$$
\eta_{t_1}-\mathcal{H}_1\ge\frac{K_2}{(1+\Lambda(y(t_1^*)))^{\beta_2}}-\mathcal{H}_1\ge\frac{K_2}{2(1+\Lambda(y(t_1^*)))^{\beta_2}}>0.
$$
Likewise, since $t_2^*\le t_2<\hat t_2$, $t_2^*\in\Gamma_2$. We then
get
$$
\xi_{t_2}-\mathcal{H}_2\ge\frac{K_1}{(1+\Lambda(x(t_2^*)))^{\beta_1}}-\mathcal{H}_2\ge\frac{K_1}{2(1+\Lambda(x(t_2^*)))^{\beta_1}}>0.
$$
By Proposition 9, for $s\le t_1$
\begin{equation}
\Lambda(x(s))\le\Lambda(x(0))e^{-2\int_0^s\eta_{\tau}-\mathcal{H}_1\mathrm{d}\tau}\le\Lambda(x(0))e^{-\frac{sK_2}{(1+\Lambda(y(t_1^*)))^{\beta_2}}};\label{12}
\end{equation}
and similarly, for $s\le t_2$,
$\Lambda(y(s))\le\Lambda(y(0))e^{-\frac{sK_1}{(1+\Lambda(x(t_2^*)))^{\beta_1}}}.$
Now we can proceed as in Proposition 5 to show neither $\hat t_1$
nor $\hat t_2$ exists. Hence, for all $t<T$,
$(1+\Lambda(y(t)))^{\beta_2}\le\frac{K_2}{2\mathcal{H}_1}=(1+\Lambda(y(0)))^{\beta_2}$,
i.e. $\Lambda(y(t))\le\Lambda(y(0))$. Combining this with expression
(\ref{12}), we obtain
$\Lambda(x(t))\le\Lambda(x(0))e^{-2t\mathcal{H}_1}$, which concludes
the first part of the proposition. The other inequality can be
derived similarly. $\Box$

\medskip
\noindent\textbf{Proof of Theorem 4}. The existence of a unique
solution for each $\omega_i\in\Omega_i$, $i=1,2$ follows from
\cite{8} (Ch. 8).

By the assumptions, we have for $\varepsilon_i>0$, $i=1,2$,
$P_i(\max_{0\le t\le T_2\vee T_3}||H_i(t)||\le \varepsilon_i)\ge
F_i(\varepsilon_i,T_2\vee T_3)$. Take
$\varepsilon_1=\mathcal{H}_1\mu$ and
$\varepsilon_2=\mathcal{H}_2\nu$ in the above expression,
respectively. By the independence of $H_1$ and $H_2$ (and let
$P:=P_1\times P_2$ as in the proof of Theorem 2), we obtain
\begin{multline}
P\big(\max_{0\le t\le T_2\vee T_3}||H_1(t)||\le \mathcal{H}_1\mu,\
\max_{0\le t\le T_2\vee T_3}||H_2(t)||\le \mathcal{H}_2\nu\big)\\
\ge F_1(\mathcal{H}_1\mu,T_2\vee T_3)F_2(\mathcal{H}_2\nu,T_2\vee
T_3).\nonumber
\end{multline}
Define $T^1(\omega_1):=\inf\{t\ge0|\ ||x(t)||\le\mu\}$ and
$T^2(\omega_2):=\inf\{t\ge0|\ ||y(t)||\le\nu\}$. It now suffices to
proof
\begin{equation}
P\big(\min\{T^1(\omega_1),T^2(\omega_2)\}\le T_2\vee T_3\big)\ge
F_1(\mathcal{H}_1\mu,T_2\vee T_3)F_2(\mathcal{H}_2\nu,T_2\vee
T_3).\label{13}
\end{equation}

Take $\omega_1\in\{\omega_1\in\Omega_1|\ \max_{0\le t\le T_2\vee
T_3}||H_1(t)||\le\mathcal{H}_1\mu\}$, and\\
$\omega_2\in\{\omega_2\in\Omega_2|\ \max_{0\le t\le T_2\vee
T_3}||H_2(t)||\le\mathcal{H}_2\nu\}$. Assume
$\min\{T^1(\omega_1),T^2(\omega_2)\}>T_2\vee T_3$. Then there are
some $\zeta>0$ such that
$\min\{T^1(\omega_1),T^2(\omega_2)\}>T_2\vee T_3+\zeta$.
Accordingly, $||x(t)||>\mu$, $||y(t)||>\nu$ on the interval
$[0,T_2\vee T_3+\zeta]$. In particular, $||x(T_2\vee T_3)||>\mu$ and
$||y(T_2\vee T_3)||>\nu$. Since $\max_{0\le t\le T_2\vee
T_3}||H_1(t)||\le\mathcal{H}_1\mu<||x(t)||\mathcal{H}_1$ and
$\max_{0\le t\le T_2\vee
T_3}||H_2(t)||\le\mathcal{H}_2\nu<||y(t)||\mathcal{H}_2$,
Proposition 10 holds for $T=T_2\vee T_3$. Then for $0\le t<T_2\vee
T_3$, we have $||x(t)||\le
||x(0)||e^{-\frac{K_2t}{2(1+\Lambda(y(0)))^{\beta_2}}}$ and
$||y(t)||\le
||y(0)||e^{-\frac{K_1t}{2(1+\Lambda(x(0)))^{\beta_1}}}$. As
$x(t),y(t)$ are continuous, we obtain $||x(T_2\vee T_3)||\le
||x(0)||e^{-\frac{K_2(T_2\vee
T_3)}{2(1+\Lambda(y(0)))^{\beta_2}}}\le\mu$ and $||y(T_2\vee
T_3)||\le ||y(0)||e^{-\frac{K_1(T_2\vee
T_3)}{2(1+\Lambda(x(0)))^{\beta_1}}}\le\nu$. Now we obtain a
contradiction, which then completes the proof of (\ref{13}). $\Box$

\medskip
\noindent{\Large\textbf{4. Conclusion}}
\smallskip

In this paper, we have studied the emergent behavior of four
dynamical systems (I(D), I(C), II(D), II(C)) in the presence of
random fluctuation contained in the environments. In all these
cases, ``nearly-emergence'' phenomena of interested objectives are
shown under certain conditions on the systems and the noises. Our
results are presented in a quite general setting and reveal some
intrinsic mechanisms of emergence which come up in a variety of
disciplines \cite{20}. We will extend the results herein onto other
dynamical systems and different kinds of random environment will be
treated in future work.

\end{document}